\pdfoutput=1

\documentclass[journal]{IEEEtran}

\usepackage{cite}
\usepackage{color}

\ifCLASSINFOpdf
  \usepackage[pdftex]{graphicx}
  \graphicspath{{../pdf/}{../jpeg/}}
\else
\fi

\usepackage{amsmath}
\usepackage{amsfonts}
\usepackage{cases}
\usepackage{amssymb}
\usepackage{bm}
\usepackage{pgfplots}
\usepackage{multirow}
\usepackage{threeparttable}

\usepackage{url}

\allowdisplaybreaks[4]

\begin{document}

\title{Robust Coordinated Transmission and Generation\\ Expansion Planning Considering Ramping Requirements
and Construction Periods}

\author{Jia~Li,~\IEEEmembership{Student~Member,~IEEE,}
        Zuyi~Li,~\IEEEmembership{Senior~Member,~IEEE,}
        Feng~Liu,~\IEEEmembership{Member,~IEEE,}
        Hongxing~Ye,~\IEEEmembership{Member,~IEEE,}
        Xuemin~Zhang,~\IEEEmembership{Member,~IEEE,}
        Shengwei~Mei,~\IEEEmembership{Fellow,~IEEE,}       
        Naichao~Chang
\thanks{J. Li, F. Liu, X. Zhang and S. Mei are with the Department of Electrical Engineering, Tsinghua University, Beijing 100084, China.}
\thanks{Z. Li and H. Ye are with Robert W. Galvin Center for Electricity Innovation at Illinois Institute of Technology, Chicago, IL 60616 USA.}
\thanks{This work was supported in part by the China State Grid Corp Science and Technology Project (SGSXDKY-DWKJ2015-001), in part by the State Key Development Program of Basic Research of China (2013CB228201).}
        }

\maketitle

\begin{abstract}
Two critical issues have arisen in transmission expansion planning with the rapid growth of wind power generation. First, severe power ramping events in daily operation due to the high variability of wind power generation pose great challenges to multi-year planning decision making. Second, the long construction periods of transmission lines may not be able to keep pace with the fast growing uncertainty due to the increasing integration of renewable energy generation. To address such issues, we propose a comprehensive robust planning model considering different resources, namely, transmission lines, generators, and FACTS devices. Various factors are taken into account, including flexibility requirement, construction period, and cost. We construct the hourly net load ramping uncertainty (HLRU) set to characterize the variation of hourly net load including wind power generation, and the annual net load duration curve uncertainty (LDCU) set for the uncertainty of normal annual net load duration curve. This results in a two-stage robust optimization model with two different types of uncertainty sets, which are decoupled into two different sets of subproblems to make the entire solution process tractable. Numerical simulations with real-world data show that the proposed model and solution method are effective to coordinate different flexible resources, rendering robust expansion planning strategies. 
\end{abstract}


\IEEEpeerreviewmaketitle

\section*{Nomenclature}

\addcontentsline{toc}{section}{Nomenclature}
\subsection{Indices}
\begin{IEEEdescription}[\IEEEusemathlabelsep\IEEEsetlabelwidth{$P^{f,\max}_{m}$/$P^{l,\max}_{ij}$}]
  \item[$h/y$] Hourly/yearly time index.
  \item[$i,j$] Bus index.
  \item[$ij$] Transmission corridor index.
  \item[$k/s$] Line/generator index.
  \item[$m/w$] FACTS device/generator type index.
  \item[$y_0$] Base year.
\end{IEEEdescription}

\subsection{Parameters}
\begin{IEEEdescription}[\IEEEusemathlabelsep\IEEEsetlabelwidth{$P^{f,\max}_{m}$/$P^{l,\max}_{ij}$}]
  \item[$a_i, b_i$] Fuel cost coefficients of an existing generator.
  \item[$a_{w}, b_{w}$] Fuel cost coefficients of a new generator.
  \item[$c_{m}/c_{w}/c_{ij}$] Investment costs of a new FACTS device/generator/line.
  \item[$d_{h}$] Time slot in an annual net load duration curve.
  \item[$n^{g,\max}_{i}$] Maximum number of new generators at each bus.
  \item[$n^{l,\max}_{ij}$/$n^{l,\min}_{ij}$] Maximum/minimum number of lines in each transmission corridor.
  \item[$P^d_{i,y,h}$] Power demand.
  \item[$P^{f,\max}_{m}$/$P^{l,\max}_{ij}$] Capacity of a FACTS device/line.
  \item[$P^{g,\max}_i/P^{gn,\max}_{w}$] Maximum power output of an existing/new generator.
  \item[$P^{g,\min}_i/P^{gn,\min}_{w}$] Minimum power output of an existing/new generator.
  \item[$R^u_{i}/R^d_{i}$] Ramp up/down limit of an existing generator.
  \item[$R^{un}_{w}/R^{dn}_{w}$] Ramp up/down limit of a new generator.
  \item[$X_{ij}$] Reactance of a transmission lines.
  \item[$y^{g}/y^{l}$] Construction period of a new generator/line.
\end{IEEEdescription}

\subsection{Variables}
\begin{IEEEdescription}[\IEEEusemathlabelsep\IEEEsetlabelwidth{$v^g_{i,y,h}/v^{gn}_{i,y,h,s,w}$}]
  \item[$n^f_{ij,y,k,m}$/$n^l_{ij,y,k}$] Binary variable indicating the installation of a new FACTS device/line.
  \item[$P^c_{ij,y,h}$] Power flow through a transmission corridor.
  \item[$P^f_{ij,y,h,k,m}$] Power injection of a FACTS device.
  \item[$P^l_{ij,y,h,k}$] Power flow in a transmission line.
  \item[$v^g_{i,y,h}/v^{gn}_{i,y,h,s,w}$] Status of an existing/new generator.
  \item[$\theta_{i,y,h}$] Phase angle.
\end{IEEEdescription}

\section{Introduction}
\IEEEPARstart{T}{he} rapid growth of wind power generation has posed new challenges to transmission expansion planning (TEP), among which two critical issues need to be addressed:  1) the uncertainty and variability of wind power require more flexibility in power system operation, particularly the ramping capability; 2) the long construction periods of transmission lines \cite{Rious2011a} may not be able to keep pace with the fast growth of wind power penetration. 

Great efforts have been devoted to the handling of wind power uncertainty in TEP. To characterize the uncertainty, stochastic optimization approaches generate a number of scenarios \cite{Carrion2007a,Maghouli2011a}, while robust optimization approaches employ uncertainty sets \cite{Jabr2013a,Jadid2013a,Kazemi2014a,ChenB.WangJ.WangL.HeY.andWang2014,Wen2015}. In order to exploit the flexibility of power system, TEP is usually coordinated with generation expansion planning (GEP) \cite{AlvarezLopez2007a,Tor2008a,Roh2009a,Pozo2013a}. In \cite{Saboori2011a}, wind power penetration is considered in a composite generation and transmission expansion model solved by a branch and bound method. A multi-objective probabilistic coordinated generation and transmission expansion framework is proposed in \cite{Heidari2015a}, and normal boundary intersection method is used to obtain the Pareto-optimal solutions. A tri-level reliability-constrained robust power system expansion planning model is proposed in \cite{Dehghan2015}, where both discrete and continuous uncertain variables are considered simultaneously. 

Note that in the state-of-the-art TEP model, operational constraints are usually formulated based on annual load duration curve (LDC), monthly load blocks \cite{Khodaei2012a}, or selected scenarios \cite{Maghouli2011a}, \cite{Garces2009a}. However, the ramping requirements have received little attention. As more and more wind power is integrated, the ramping capability of a power system will become a key limiting factor to its capability of accommodating variable wind generation. This motivated us to explicitly incorporate ramping requirements into the TEP model in this paper. For the convenience of discussion, we consider only the net load which is defined as the remaining system load not served by the wind power generation \cite{Lannoye2012e}. 

New technologies in transmission network, such as flexible AC transmission systems (FACTS), have also been considered in TEP. Phase-shifter transformer is considered as an element in TEP to extend the utilization of classical components in \cite{Miasaki2012a}. An investment valuation approach is proposed in real option analysis framework to assess the option value of FACTS in TEP \cite{Blanco2011a}. In \cite{Konstantelos2015a}, energy storage, demand-side management, and phase-shifting transformer are incorporated into TEP based on a stochastic framework to investigate the potential of these non-conventional assets in accommodating renewable energy. It is worth noting that FACTS can not only provide extra flexibility to serve as a supplement in TEP decision making, but also enable new opportunities to coordinate the construction process, due to their relatively short construction periods, as we will reveal in this paper.

Aiming at addressing the two issues mentioned above, we propose a comprehensive robust planning model incorporating different flexible resources and considering different uncertainties. The main contributions of this paper are threefold:

1) Proposing a comprehensive multi-year planning model, which incorporates three typical resources of flexibility: transmission lines, generators, and FACTS devices. Construction periods are considered to investigate the impacts of FACTS devices. The planning strategy is based on an overall consideration of various factors, including flexibility requirement, construction period, and cost. Numerical results reveal that FACTS devices can help cope with uncertainty and coordinate resources with different construction periods, providing insights to power system planners about the coordination of different resources.

2) Modeling the uncertainty of ramping in the multi-year planning model, which is a trade-off between accuracy and tractability. With increasing integration of renewable energy generation, the operational flexibility of the system may not be sufficient due to the lack of ramping capability. In the literature, the uncertainty of annual net load duration curve (LDCU) is usually taken into account using typical scenarios or uncertainty sets, while scenarios have been used to model the uncertain hourly load variation as well as the ramping effect \cite{Jabr2014a}. However, in our multi-year planning problem, a few scenarios are not enough to characterize the uncertainty of ramping, since the precision of hourly load curve forecasts several years ahead is low currently. Meanwhile, a large number of scenarios will make the optimization problem intractable. Therefore, an additional uncertainty set is introduced to describe the hourly net load ramping uncertainty (HLRU) without causing computational intractability. Real load data are used to verify the validity of the proposed model. The results show that the proposed model renders a reliable planning strategy while avoiding computational intractability.

3) Improving computational efficiency of the two-stage robust optimization problem with two coupled uncertainties. Even though applying uncertainty sets renders a tractable mixed-integer linear programming (MILP) problem, the computational burden is still high due to the consideration of the two types of uncertainties, especially for large-scale systems. In this paper, the two types of uncertainties are decoupled into two subproblems to improve the efficiency of the Column-and-Constraint Generation (C\&CG) \cite{Zeng2013} method in our case. The two subproblems are solved using different methods, and the Relax-and-Enforce Decomposition (RED) \cite{Ye2015c} technique is applied to make a temporal decomposition to further reduce computational burden. Tests on the IEEE 118-bus system and a real-world system show that the RED technique may accelerate the solution process as much as one order of magnitude compared with the standard C\&CG method without a temporal decomposition.

The rest of the paper is organized as follows. Section \ref{sec_formulation} gives the comprehensive TEP model. The solution approach is presented in Section \ref{sec_solution}. In Section \ref{sec_case}, case studies using real-world data are conducted. Section \ref{sec_conclusion} draws the conclusions.

\section{Problem Formulation} \label{sec_formulation}
The aim of the proposed robust planning model is to find the least total cost, including the investment cost and the operation cost of the base-case scenario, over two feasibility sets defined by uncertainties, one of which refers to LDCU, denoted by $\mathcal{F}^d$, and the other refers to HLRU, denoted by $\mathcal{F}^r$. The base-case scenario is built according to the forecasted data, while generators can be re-dispatched and FACTS devices can be adjusted to control the power flow when uncertainty is revealed. Meanwhile, the re-dispatch cost for accommodating deviations from the base-case scenario, known as recourse cost \cite{Ye2015c}, should be limited under an acceptable level. Note that in this paper the planning strategy is assumed always feasible without load shedding, which may be stronger than the traditional TEP model where load shedding is allowed. The conservativeness of the proposed model, however, can be controlled by adjusting the recourse cost and uncertainty sets. Besides, one can further consider the constraints on the acceptable amount of load shedding in the our model, which will be presented in our future work.

\subsection{Uncertainty Modeling and Decoupling}
\begin{figure}[htb]
\centering
\includegraphics[width=3.6in]{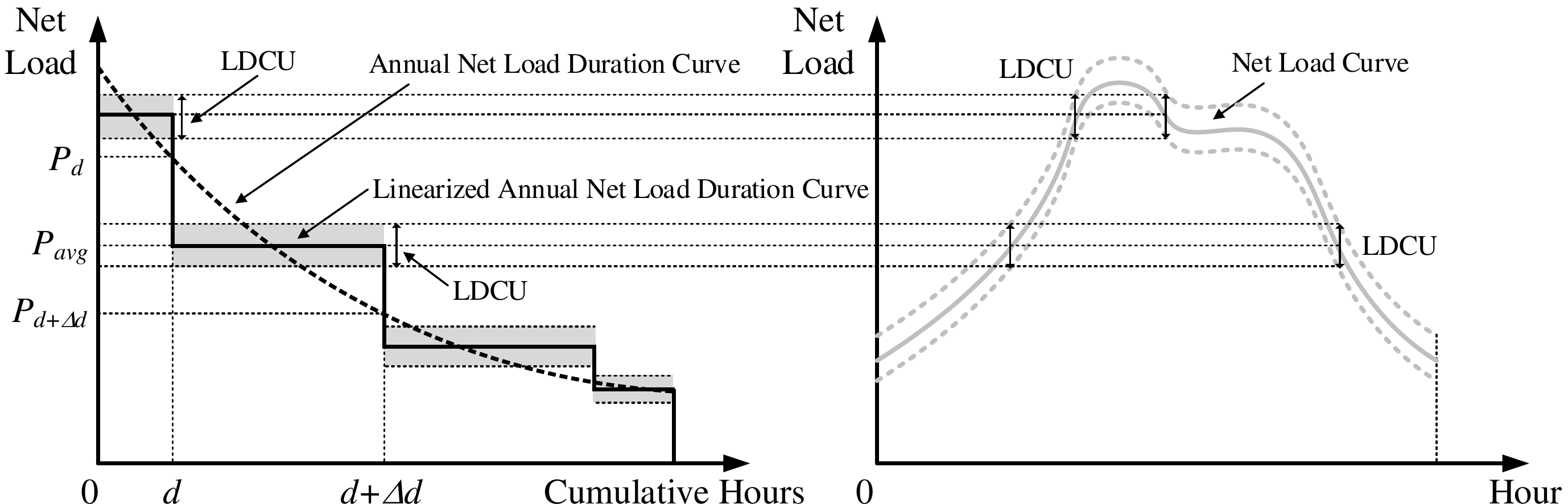}
\caption{Load duration curve uncertainty (LDCU).}
\label{fig_ldcu}
\end{figure}

As mentioned previously, we depict the uncertainty directly based on net load for the convenience of discussion. Fig.\ref{fig_ldcu} illustrates the LDCU. The annual net load duration curve is linearized by dividing the whole time period into several time slots and averaging the net load within each time slot. For example, there are $d$ hours during a year when the net load level is higher than $P_d$, and $d+\Delta d$ hours when the net load level is higher than $P_{d+\Delta d}$. Assume that $P^{avg}_d$ is the average net load of all the net loads within the time slot $[d,d+\Delta d]$ and that $P^{avg}_d$ lasts for $\Delta d$ hours during a year. Then, the LDCU indicates that each average net load level $P^{avg}_d$ may be higher or lower than the forecasted value due to unavoidable forecasting error of the associated net load duration curve.

\begin{figure}[htb]
\centering
\includegraphics[width=3.6in]{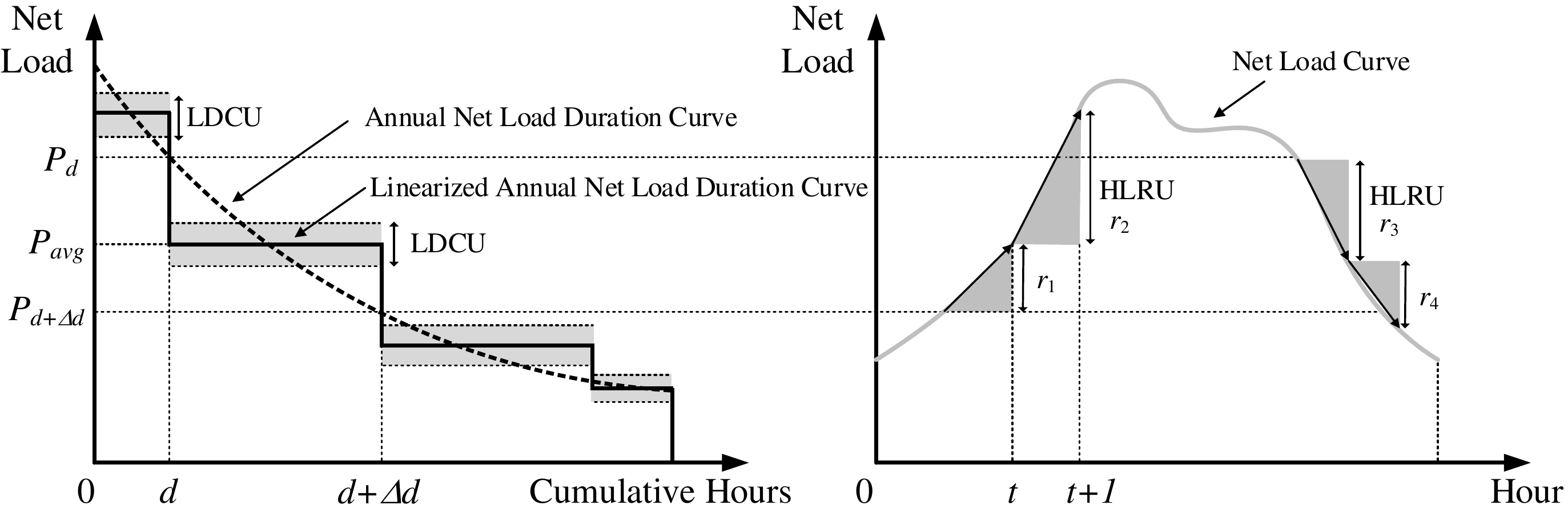}
\caption{Two types of uncertainties (LDCU and HLRU).}
\label{fig_hlru}
\end{figure}

Fig.\ref{fig_hlru} illustrates the relationship between the two types of uncertainty. It should be noted that the annual net load duration curve only arranges all the net load levels in a descending order of the magnitude, regardless of the temporal variation between two time slots. As a result, the hourly variation of  net load  is ignored. The traditional method to consider such variation is to simply use some typical daily load curves, as simulating all possible daily load curves is intractable in a multi-year planning model. However, such treatment may be insufficient to cover all the variations when a large amount of wind power is integrated with large uncertainty. In this context, the HLRU is introduced as a trade-off to describe the hourly variation of net load, hence bridging the gap between the annual net load duration curve and the daily net load curve. 

Note that the two types of uncertainty are coupled to each other, creating considerable difficulty in problem solving. In fact, 8760 hours net load are sorted to make an annual net load duration curve, indicating that each net load level in the annual net load duration curve corresponds to a net load level at a specific hour in a certain day. Suppose that the specific hour is $t$. Then, the LDCU represents the uncertainty of net load at hour $t$, while the HLRU indicates the net load variation from hour $t$ to $t+1$. Note that there are infinite combinations of net load level at hour $t$ and ramping from hour $t$ to $t+1$, when both the LDCU and the HLRU are taken into account. Hence, we decouple the LDCU from the HLRU and assume that the ramping approximately starts from each average net load level $P^{avg}_d$. Then we examine the annual net load curve to record those hourly ramping values where the ramping events start from the net load levels within the range $[P_{d+\Delta d},P_d]$. For example, ramping value $r_1$, $r_2$, $r_3$, $r_4$ are recorded in Fig. \ref{fig_hlru}. Last, we select the highest ramping up/down value as the upper/lower bound of the HLRU. In this way, the two types of uncertainties are decoupled, resulting in two uncertainty sets, $\mathcal{U}^d$ and $\mathcal{U}^r$, as follows.
\begin{equation}\label{eq_ul}
  \begin{aligned}
    \mathcal{U}^d & := \{ \epsilon_{i,y,h}^d\in\mathbb{R}: \forall i,y,h, \\
    & |\epsilon_{i,y,h}^d| \le \bar{u}_{i,y,h}^d, \sum_i\frac{|\epsilon_{i,y,h}^d|}{\bar{u}_{i,y,h}^d}\le\Lambda_{y,h} \}
  \end{aligned}
\end{equation}
\begin{equation}\label{eq_ur}
  \begin{aligned}
    \mathcal{U}^r & := \{ \epsilon_{i,y,h}^r\in\mathbb{R}: \forall i,y,h, \underline{u}_{i,y,h}^r\le\epsilon_{i,y,h}^r\le\bar{u}_{i,y,h}^r \}
  \end{aligned}
\end{equation}

The characteristics of these two uncertainty sets and the resources available to accommodate the uncertainty are different in the following aspects: 1) the absolute values of the upper and lower bounds of $\epsilon_{i,y,h}^d$ are assumed to be the same, indicating the equal possibility of over-forecast and under-forecast, while the upper limit $\bar{u}_{i,y,h}^r$ and lower limit $\underline{u}_{i,y,h}^r$ of $\epsilon_{i,y,h}^r$ , which are related to the maximum ramp-up and ramp-down events, may be different; 2) the budget of uncertainty $\Lambda_{y,h}$ \cite{Bertsimas2013} is employed in (\ref{eq_ul}) to control the conservativeness related to forecasting error; 3) the intra-hour flexible resources are required to accommodate the LDCU, such as 10-minute spinning reserve, while the inter-hour flexible resources are applied to address the HLRU, such as hourly ramping capability of generators.

\begin{figure}[htb]
\centering
\includegraphics[width=3in]{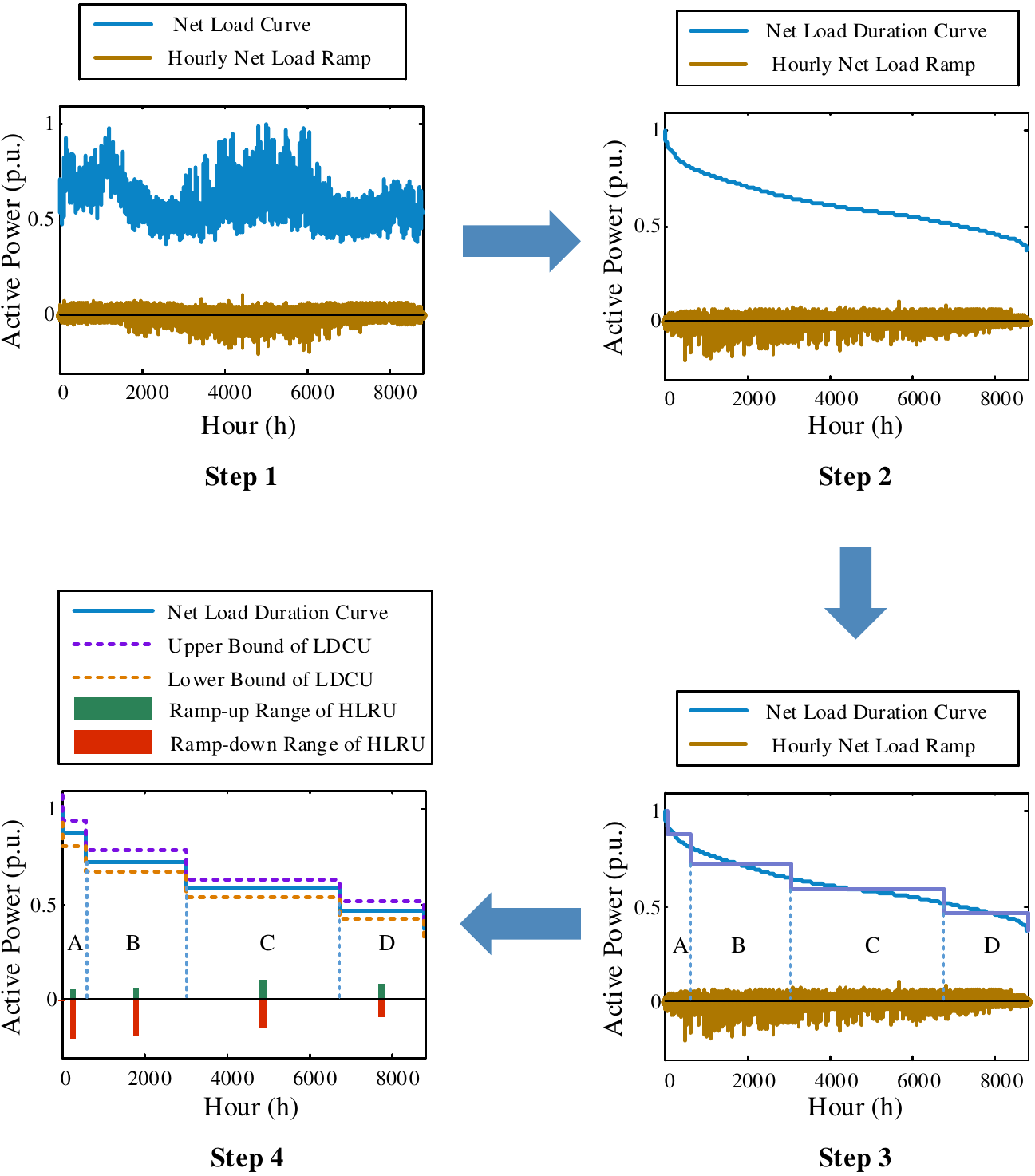}
\caption{Process of constructing the LDCU and HLRU sets.}
\label{fig_example}
\end{figure}

Figure. \ref{fig_example} depicts an illustrative example of the process of constructing the LDCU and HLRU sets based on an annual load curve. The process can be divided into the following four steps:

\textbf{Step 1: Constructing hourly net load ramping events from an annual net load curve.} The annual net load curve contains 8760 hours net loads. Hourly net load ramping event is the difference of the net loads between two adjacent hours. Each hourly net load ramping event is attached to the hourly net load level where the ramping starts from.

\textbf{Step 2: Constructing an annual net load duration curve and arranging the corresponding hourly net load ramping events.} The annual net load duration curve is derived by sorting all the net loads in a descending order. Accordingly, the hourly net load ramping events, which are attached to net load levels in Step 1, are arranged.

\textbf{Step 3: Linearizing the annual net load duration curve.} The annual net load duration curve is divided into several slots. In this example, it is divided into 4 slots, i.e., A, B, C, D. Then the average value of the net load levels within each slot is used to linearize the annual net load duration curve as a stepwise function. Accordingly, all the ramping events are divided into the same 4 slots, i.e., A, B, C, D. For instance, in slot A, suppose there are 100 different net load levels. For linearizion, the average value of these 100 net load levels represents the net load level of slot A. And the hourly net load ramping events attached to these 100 net load levels belong to slot A.

\textbf{Step 4: Constructing uncertainty sets.} Based on the average net load of each slot, the upper and lower bounds of LDCU set are constructed according to the forecast error of net load curves. Meanwhile, the range of HLRU set for each slot is derived from the maximum hourly net load ramp-up and ramp-down events within each slot. For instance, in slot A, the 100 hourly ramping events are examined to select the maximum ramp-up value and the maximum ramp-down value as the upper bound and the lower bound of the HLRU set associated with slot A, respectively.

This method bridges the gap between the long-term uncertainty of planning and the short-term uncertainty of operation, and is also a trade-off between accuracy and tractability.

\subsection{Objective Function}
Let $\boldsymbol{n}^f$, $\boldsymbol{n}^g$, $\boldsymbol{n}^l$ denote the vectors of binary variables indicating the installation of new FACTS devices, generators, and lines, respectively; let vector $\boldsymbol{v}$ denote the binary variables indicating the status of generator; let vectors $\boldsymbol{p}^c$, $\boldsymbol{p}^f$, $\boldsymbol{p}^g$, $\boldsymbol{p}^l$ denote active power associated with transmission corridor, FACTS, generator and line, respectively; let vector $\boldsymbol{\theta}$ denote the phase angle; let $\boldsymbol{\epsilon}^d$ and $\boldsymbol{\epsilon}^r$ be uncertainty variables associated with LDCU and HLRU, respectively. The objective (\ref{eq_obj}) is to minimize the total investment and operation cost in the base-case scenario shown below.
\begin{equation}
  \min_{(\boldsymbol{n}^f,\boldsymbol{n}^g,\boldsymbol{n}^l,\boldsymbol{v},\boldsymbol{p}^g)\in\mathcal{F}^d\cap\mathcal{F}^r} \boldsymbol{a}^{\top}\boldsymbol{n}^l + \boldsymbol{b}^{\top}\boldsymbol{n}^f + \boldsymbol{c}^{\top}\boldsymbol{n}^g + \boldsymbol{d}^{\top}\boldsymbol{v} + \boldsymbol{e}^{\top}\boldsymbol{p}^g \label{eq_obj}
\end{equation}
where
\begin{align}
  \boldsymbol{a}^{\top}\boldsymbol{n}^l &= \sum_{y}\sum_{ij}\sum_{k}\frac{c_{ij}(n_{ij,y,k}^l - n_{ij,y-1,k}^l)}{(1+D)^{(y-y_0 - y^l)}} \label{eq_cost_line} \\
  \boldsymbol{b}^{\top}\boldsymbol{n}^f &= \sum_{y}\sum_{ij}\sum_{k}\sum_{m}\frac{c_{m}(n_{ij,y,k,m}^f - n_{ij,y-1,k,m}^f)}{(1+D)^{(y-y_0)}} \label{eq_cost_facts} \\
  \boldsymbol{c}^{\top}\boldsymbol{n}^g &= \sum_{y}\sum_{i}\sum_{s}\sum_{w}\frac{c_{w}(n_{i,y,s,w}^g - n_{i,y-1,s,w}^g)}{(1+D)^{(y - y_0 - y^g)}} \label{eq_cost_gen} \\
  \boldsymbol{d}^{\top}\boldsymbol{v} + \boldsymbol{e}^{\top}\boldsymbol{p}^g &= \sum_{y}\sum_{h}\sum_{i}\left\{\frac{d_{h}(a_i v_{i,y,h} + b_i P^g_{i,y,h})}{(1+D)^{(y - y_0)}} \right. \nonumber \\
  & \left.+ \sum_{s}\sum_{w}\frac{d_{h}(a_{w} v^{gn}_{i,y,h,s,w} + b_{w} P^{gn}_{i,y,h,s,w})}{(1+D)^{(y - y_0 - y^g)}} \right\}. \label{eq_cost_op}
\end{align}
Assume that a new line takes $y^l$ years to build and the construction is finished in year $t$. Then, the investment of the line has to be made in year $t-y^l$. Accordingly, $t-y_0-y^l$ is the period with a discount rate $D$, as shown in (\ref{eq_cost_line}). If $y^l$ equals 0, it means that the new line is put into use in the same year when the investment is made. The investment cost of a generator is formulated similarly (\ref{eq_cost_gen}). It is assumed that a FACTS device is installed in the same year when it is invested (\ref{eq_cost_facts}). The operation cost of the base-case scenario is formulated using linear production cost function as in (\ref{eq_cost_op}).

\subsection{Incorporating FACTS into TEP} \label{sec_facts}
We incorporate FACTS devices into the TEP model using the power injection model \cite{WangX.2002}. The nodal power balance equation is formulated as (\ref{eq_pb}), where a FACTS device is formulated as two power injections, which have the same amount, but opposite signs, located at each bus of a transmission corridor.
\begin{equation}
  \begin{aligned}
    \sum_jP^c_{ij,y,h} & = P^g_{i,y,h} + \sum_{s}\sum_{w}P^{gn}_{i,y,h,s,w} - P^d_{i,y,h} \\
    & + \sum_{j,i<j}\sum_{k}\sum_{m}P^f_{ij,y,h,k,m} - \sum_{j,i>j}\sum_{k}\sum_{m}P^f_{ji,y,h,k,m} \label{eq_pb}
  \end{aligned}
\end{equation}
The power injection of FACTS is constrained by (\ref{eq_p_facts})
\begin{equation}
        |P^f_{ij,y,h,k,m}| \le P^{f,\max}_{m}n^f_{ij,y,k,m}. \label{eq_p_facts}
\end{equation}
The transmission capacity is constrained by
\begin{equation}
        |P^l_{ij,y,h,k} - \sum_m P^f_{ij,y,h,k,m}| \le P^{l,\max}_{ij}n^l_{ij,y,k}.
\end{equation}
The DC power flow is enforced by
\begin{equation}
    |\theta_{i,y,h} - \theta_{j,y,h} - X_{ij}P^l_{ij,y,h,k}| \le 2\theta^{\max}(1 - n^l_{ij,y,k}) \label{eq_pf}
\end{equation}
where the phase angles satisfy
\begin{equation}
  |\theta_{i,y,h}| \le \theta^{\max} \label{eq_theta}
\end{equation}
and the maximum phase angle $\theta^{\max}$ is set as $\pi/2$ \cite{Vinasco2011a}.
Constraint (\ref{eq_nf_nl}) allows the installation of a FACTS device only in the existing line, but no more than one FACTS device is allowed in each line.
\begin{equation}
  \sum_{k}\sum_{m}n^f_{ij,y,k,m} \le \sum_{k}n^l_{ij,y,k} \label{eq_nf_nl}
\end{equation}
Once a FACTS device is installed, it will exist during the rest of the planning horizon (\ref{eq_nf_bo}).
\begin{equation}
  n^f_{ij,y,k,m} \ge n^f_{ij,y-1,k,m} \label{eq_nf_bo}
\end{equation}
Besides, new FACTS devices in each corridor are installed sequentially (\ref{eq_nf_seq}).
\begin{equation}
  n^f_{ij,y,k,m} \le n^f_{ij,y,k-1,m} \label{eq_nf_seq}
\end{equation}
At last, only one type of FACTS device is allowed to be installed in each line (\ref{eq_nf_type}).
\begin{equation}
  \sum_{m}n^f_{ij,y,k,m} \le 1 \label{eq_nf_type}
\end{equation}

\subsection{Other Constraints in Base-case Scenario}
Each transmission corridor has a minimum limit and a maximum limit to the number of lines (\ref{eq_nl_minmax}). The minimum number refers to the number of existing lines in each corridor. Constraint (\ref{eq_nl_t0}) makes sure no new line is built before the first construction period.
\begin{equation}
  n_{ij}^{l,\min} \le \sum_{k}n_{ij,y,k}^l \le n_{ij}^{l,\max} \label{eq_nl_minmax}
\end{equation}
\begin{equation}
  \sum_{k}n_{ij,y,k}^l \le n_{ij}^{l,\min}, \forall y < y_0 + y^l \label{eq_nl_t0}
\end{equation}
Constraint (\ref{eq_nl_bo}) indicates that once a new line is built, it will exist during the planning horizon.
\begin{equation}
  n^l_{ij,y,k} \ge n^l_{ij,y-1,k} \label{eq_nl_bo}
\end{equation}
Besides, the new lines in each transmission corridor are built sequentially (\ref{eq_nl_seq}).
\begin{equation}
  n^l_{ij,y,k} \le n^l_{ij,y,k-1} \label{eq_nl_seq}
\end{equation}
The total power flow through a transmission corridor is the sum of the power flows on all lines in that corridor (\ref{eq_p_cor}).
\begin{equation}
  P_{ij,y,h}^c = \sum_{k}P_{ij,y,h,k}^l \label{eq_p_cor}
\end{equation}
The maximum number of new generators is limited by (\ref{eq_ng_max}).
\begin{equation}
  \sum_{s}\sum_{w}n^g_{i,y,s,w} \le n^{g,\max}_i \label{eq_ng_max}
\end{equation}
Once a new generator is installed, it will exist during the planning horizon (\ref{eq_ng_bo}). 
\begin{equation}
  n^g_{i,y,s,w} \ge n^g_{i,y-1,s,w} \label{eq_ng_bo}
\end{equation}
Besides, the new generators at each bus are built sequentially (\ref{eq_ng_seq}), and only one type of new generator would be installed each time (\ref{eq_ng_type}).
\begin{equation}
  n^g_{i,y,s,w} \le n^g_{i,y,s-1,w} \label{eq_ng_seq}
\end{equation}
\begin{equation}
  \sum_{w}n^g_{i,y,s,w} \le 1 \label{eq_ng_type}
\end{equation}
Constraint (\ref{eq_ng_v}) indicates that a new generator can only be dispatched after installation (\ref{eq_ng_v}).
\begin{equation}
  v^{gn}_{i,y,h,s,w} \le n^g_{i,y,s,w} \label{eq_ng_v}
\end{equation}
Constraints (\ref{eq_p_g}) and (\ref{eq_p_gn}) enforce the capacity limits of the existing generators and new generators, respectively.
\begin{align}
  & P^{g,\min}_i v_{i,y,h} \le P^g_{i,y,h} \le P^{g,\max}_i v_{i,y,h} \label{eq_p_g} \\
  & P^{gn,\min}_{w} v^{gn}_{i,y,h,s,w} \le P^{gn}_{i,y,h,s,w} \le P^{gn,\max}_{w} v^{gn}_{i,y,h,s,w} \label{eq_p_gn}
\end{align}

\subsection{Robust Planning Model}
The robust TEP model is formulated in a compact form below, where bold symbols except for the variables mentioned before are constant matrices or vectors. Vectors with 1 and 2 in subscripts represent variables associated with LDCU and HLRU, respectively.
\begin{align}
  \min_{(\boldsymbol{n}^f,\boldsymbol{n}^g,\boldsymbol{n}^l,\boldsymbol{v},\boldsymbol{p}^g)\in\mathcal{F}^d\cap\mathcal{F}^r} & \boldsymbol{a}^{\top}\boldsymbol{n}^l + \boldsymbol{b}^{\top}\boldsymbol{n}^f + \boldsymbol{c}^{\top}\boldsymbol{n}^g + \boldsymbol{d}^{\top}\boldsymbol{v} + \boldsymbol{e}^{\top}\boldsymbol{p}^g \label{eq_compact_obj} \\
  \text{s.t.}\quad & \boldsymbol{A}^l\boldsymbol{n}^l \le \boldsymbol{f} \label{eq_compact_nl} \\
  & \boldsymbol{B}^c\boldsymbol{p}^c + \boldsymbol{C}^l\boldsymbol{p}^l = \boldsymbol{0} \label{eq_compact_p_cor} \\
  & \boldsymbol{D}^l\boldsymbol{n}^l + \boldsymbol{E}^l\boldsymbol{p}^l + \boldsymbol{Lp}^f \le \boldsymbol{0} \label{eq_compact_p_l} \\
  & \boldsymbol{F}^p\boldsymbol{p}^c + \boldsymbol{G}^p\boldsymbol{p}^g + \boldsymbol{H}^p\boldsymbol{p}^f = \boldsymbol{g} \label{eq_compact_facts_pb} \\
  & \boldsymbol{D}^f\boldsymbol{n}^f + \boldsymbol{E}^f\boldsymbol{p}^f \le \boldsymbol{0} \label{eq_compact_facts_cap}\\
  & \boldsymbol{Jn}^l + \boldsymbol{Kp}^l + \boldsymbol{M\theta} \le \boldsymbol{h} \label{eq_compact_facts_pf} \\
  & \boldsymbol{A}^{\theta}\boldsymbol{\theta} \le \boldsymbol{k} \label{eq_compact_facts_theta} \\
  & \boldsymbol{N}^l\boldsymbol{n}^l + \boldsymbol{Q}^f\boldsymbol{n}^f \le \boldsymbol{0} \label{eq_compact_facts_nl_nf}\\
  & \boldsymbol{A}^f\boldsymbol{n}^f \le \boldsymbol{l} \label{eq_compact_facts_nf} \\
  & \boldsymbol{A}^g\boldsymbol{n}^g \le \boldsymbol{m} \label{eq_compact_ng} \\
  & \boldsymbol{N}^g\boldsymbol{n}^g + \boldsymbol{Q}^v\boldsymbol{v} \le \boldsymbol{0} \label{eq_compact_ng_v} \\
  & \boldsymbol{D}^v\boldsymbol{v} + \boldsymbol{E}^v\boldsymbol{p}^g \le \boldsymbol{0} \label{eq_compact_uc}
\end{align}
where
\begin{align}
  \mathcal{F}^d := \{ & (\boldsymbol{n}^f,\boldsymbol{n}^g,\boldsymbol{n}^l,\boldsymbol{v},\boldsymbol{p}^g):\forall\boldsymbol{\epsilon}^d\in\mathcal{U}^d, \exists\boldsymbol{p}^c_1,\boldsymbol{p}^f_1,\boldsymbol{p}^g_1,\boldsymbol{p}^l_1,\boldsymbol{\theta}_1  \nonumber \\
  \text{s.t.}\quad & \boldsymbol{B}^c\boldsymbol{p}^c_1 + \boldsymbol{C}^l\boldsymbol{p}^l_1 = \boldsymbol{0} \\
  & \boldsymbol{D}^l\boldsymbol{n}^l + \boldsymbol{E}^l\boldsymbol{p}^l_1 + \boldsymbol{Lp}^f_1 \le \boldsymbol{0} \\
  & \boldsymbol{F}^p\boldsymbol{p}^c_1 + \boldsymbol{G}^p\boldsymbol{p}^g_1 + \boldsymbol{H}^p\boldsymbol{p}^f_1 + \boldsymbol{R}\boldsymbol{\epsilon}^d = \boldsymbol{g} \label{eq_compact_spd_pb} \\
  & \boldsymbol{D}^f\boldsymbol{n}^f + \boldsymbol{E}^f\boldsymbol{p}^f_1 \le \boldsymbol{0} \\
  & \boldsymbol{Jn}^l + \boldsymbol{Kp}^l_1 + \boldsymbol{M\theta}_1 \le \boldsymbol{h} \\
  & \boldsymbol{A}^{\theta}\boldsymbol{\theta}_1 \le \boldsymbol{k} \\
  & \boldsymbol{D}^v\boldsymbol{v} + \boldsymbol{E}^v\boldsymbol{p}^g_1 \le \boldsymbol{0} \\
  & \boldsymbol{e}^{\top}(\boldsymbol{p}^g_1 - \boldsymbol{p}^g) \le c^d\}, \label{eq_compact_spd_rc}
\end{align}
\begin{align}
  \mathcal{F}^r := \{ & (\boldsymbol{n}^f,\boldsymbol{n}^g,\boldsymbol{n}^l,\boldsymbol{v},\boldsymbol{p}^g):\forall\boldsymbol{\epsilon}^r\in\mathcal{U}^r, \exists\boldsymbol{p}^c_2,\boldsymbol{p}^f_2,\boldsymbol{p}^g_2,\boldsymbol{p}^l_2,\boldsymbol{\theta}_2  \nonumber \\
  \text{s.t.}\quad & \boldsymbol{B}^c\boldsymbol{p}^c_2 + \boldsymbol{C}^l\boldsymbol{p}^l_2 = \boldsymbol{0} \\
  & \boldsymbol{D}^l\boldsymbol{n}^l + \boldsymbol{E}^l\boldsymbol{p}^l_2 + \boldsymbol{Lp}^f_2 \le \boldsymbol{0} \\
  & \boldsymbol{F}^p\boldsymbol{p}^c_2 + \boldsymbol{G}^p\boldsymbol{p}^g_2 + \boldsymbol{H}^p\boldsymbol{p}^f_2 + \boldsymbol{R}\boldsymbol{\epsilon}^r = \boldsymbol{g} \label{eq_compact_spr_pb} \\
  & \boldsymbol{D}^f\boldsymbol{n}^f + \boldsymbol{E}^f\boldsymbol{p}^f_2 \le \boldsymbol{0} \\
  & \boldsymbol{Jn}^l + \boldsymbol{Kp}^l_2 + \boldsymbol{M\theta}_2 \le \boldsymbol{h} \\
  & \boldsymbol{A}^{\theta}\boldsymbol{\theta}_2 \le \boldsymbol{k} \\
  & \boldsymbol{D}^v\boldsymbol{v} + \boldsymbol{E}^v\boldsymbol{p}^g_2 \le \boldsymbol{0} \\
  & \boldsymbol{S}\boldsymbol{p}^g + \boldsymbol{T}\boldsymbol{p}^g_2 \le \boldsymbol{\Delta} \}. \label{eq_compact_spr_ramp}
\end{align}

Constraint (\ref{eq_compact_nl}) refers to the limits of new lines (\ref{eq_nl_minmax})-(\ref{eq_nl_seq}). Equation (\ref{eq_compact_p_cor}) denotes (\ref{eq_p_cor}). Constraints (\ref{eq_compact_p_l})-(\ref{eq_compact_facts_nf}) represent the constraints corresponding to FACTS devices, which are explained in detail in Section \ref{sec_facts}. Constraint (\ref{eq_compact_ng}) is associated with the limits of new generators (\ref{eq_ng_max})-(\ref{eq_ng_type}). Constraint (\ref{eq_compact_ng_v}) is the compact form of (\ref{eq_ng_v}). Constraint (\ref{eq_compact_uc}) refers to (\ref{eq_p_g})-(\ref{eq_p_gn}).

The solution to (\ref{eq_compact_nl})-(\ref{eq_compact_uc}) in the base-case scenario should be feasible for both feasibility sets, i.e. $\mathcal{F}^d$ and $\mathcal{F}^r$, with flexible resources accommodating uncertainty. When $\boldsymbol{\epsilon}^d$ varies within the uncertainty set $\mathcal{U}^d$, flexible resources are re-dispatched while power balance constraint (\ref{eq_compact_spd_pb}) is satisfied, and the recourse cost is within an acceptable level $c^d$ (\ref{eq_compact_spd_rc}), which is reformulated in detail as (\ref{eq_spd_rc}). 
\begin{align}
  & \sum_{y}\sum_{h}\sum_{i}\sum_{s}\sum_{w}\frac{d_{h}b_{w}(P^{gn}_{1,i,y,h,s,w} - P^{gn}_{i,y,h,s,w})}{(1+D)^{(y - y_0 - y^g)}} \nonumber \\
  & + \sum_{y}\sum_{h}\sum_{i}\frac{d_{h}b_i(P^g_{1,i,y,h} - P^g_{i,y,h})}{(1+D)^{(y - y_0)}} \le c^d \label{eq_spd_rc}
\end{align}
Note that there is no ramping constraint in $\mathcal{F}^d$, since this set is to check if the planning strategy and unit commitment schedule are feasible within a certain recourse cost limit, while the economic dispatch is performed based on the plausible load levels. $\boldsymbol{\epsilon}^r$ is introduced into the power balance equation in $\mathcal{F}^r$ (\ref{eq_compact_spr_pb}). When $\boldsymbol{\epsilon}^r$ varies within uncertainty set $\mathcal{U}^r$, the generators are re-dispatch within ramping limits (\ref{eq_compact_spr_ramp}), which are reformulated in detail as below. 
\begin{align}
  P^g_{2,i,y,h} - P^g_{i,y,h} & \le R^u_{i} \label{eq_spr_ru} \\
  P^g_{2,i,y,h} - P^g_{i,y,h} & \ge -R^d_{i} \label{eq_spr_rd} \\
  P^{gn}_{2,i,y,h,s,w} - P^{gn}_{i,y,h,s,w} & \le R^{un}_{w} \label{eq_spr_run} \\
  P^{gn}_{2,i,y,h,s,w} - P^{gn}_{i,y,h,s,w} & \ge -R^{dn}_{w} \label{eq_spr_rdn}
\end{align} 
The recourse cost constraint is omitted in $\mathcal{F}^r$ since this set is to ensure that all the possible net load ramping events can be accommodated without any physical limit violation. Other constraints in $\mathcal{F}^d$ and $\mathcal{F}^r$ are similar to those in the base-case scenario.

\section{Solution Methodology} \label{sec_solution}
The solution method is based on the C\&CG framework\cite{Zeng2013}. First, by decoupling the two types of uncertainties, the proposed model is formulated as a master problem and two subproblems. Second, taking advantage of the dual theory and the extreme point formulation, the subproblems can be reformulated as mixed-integer linear programming (MILP) problems. Third, to improve computational efficiency, the RED approach is applied to decompose the subproblem in a temporal manner.

\subsection{Column-and-Constraint Generation Method (C\&CG)}
By applying C\&CG method, the proposed model is reformulated a master problem and two subproblems. The master problem (MP) is defined in a compact form as below.
\begin{align}
  \text{(MP)}& \quad\min_{\boldsymbol{x},\boldsymbol{p}^g}\quad \boldsymbol{r}^{\top}\boldsymbol{x}+\boldsymbol{e}^{\top}\boldsymbol{p}^g \\
  \text{s.t.}\quad & \boldsymbol{Ax}+\boldsymbol{Bp}+\boldsymbol{C\theta}\le\boldsymbol{s} \\
  & \boldsymbol{D}\boldsymbol{x} + \boldsymbol{E}\boldsymbol{p}^{\tau} + \boldsymbol{F}\boldsymbol{\theta}^{\tau} + \boldsymbol{G}\boldsymbol{\epsilon}^{d,\tau} \le \boldsymbol{w}, \forall \tau\in\mathcal{T}  \\
  & \boldsymbol{e}^{\top}(\boldsymbol{p}^{g,\kappa} - \boldsymbol{p}^g) \le c^d, \forall \kappa\in\mathcal{K}  \\
  & \boldsymbol{D}\boldsymbol{x} + \boldsymbol{E}\boldsymbol{p}^{\kappa} + \boldsymbol{F}\boldsymbol{\theta}^{\kappa} + \boldsymbol{G}\boldsymbol{\epsilon}^{r,\kappa} \le \boldsymbol{w}, \forall \kappa\in\mathcal{K}   \\
  & \boldsymbol{S}\boldsymbol{p}^g + \boldsymbol{T}\boldsymbol{p}^{g,\kappa} \le \boldsymbol{\Delta}, \forall \kappa\in\mathcal{K}  
\end{align}
where $\boldsymbol{x}$ denotes all the binary variables, while $\boldsymbol{p}$ and $\boldsymbol{\theta}$ represent variables associated with active power and phase angles, respectively. $\mathcal{T}$, $\mathcal{K}$ are the index sets for uncertainty points $\boldsymbol{\epsilon}^{d,\tau}$ and $\boldsymbol{\epsilon}^{r,\kappa}$, respectively, which are generated in subproblems during iterations. Variable $\boldsymbol{p}^{\tau}$, $\boldsymbol{p}^{g,\tau}$ and $\boldsymbol{\theta}^{\tau}$ are associated with $\boldsymbol{\epsilon}^{d,\tau}$, while $\boldsymbol{p}^{\kappa}$, $\boldsymbol{p}^{g,\kappa}$ and $\boldsymbol{\theta}^{\kappa}$ are associated with $\boldsymbol{\epsilon}^{r,\kappa}$. The master problem is an MILP problem that can be solved efficiently by using commercial solvers.

The subproblem related to LDCU (SPD) is formulated as follows with subscription 1 indicating associated variables.
\begin{equation}
  \text{(SPD)}\quad\max_{\boldsymbol{\epsilon}^d\in\mathcal{U}^d}\min_{(\boldsymbol{s}_1^{+},\boldsymbol{s}_1^{-},\boldsymbol{p}_1)\in\mathcal{D}(\boldsymbol{\epsilon}^d)}\boldsymbol{1}^{\top}\boldsymbol{s}_1^{+} + \boldsymbol{1}^{\top}\boldsymbol{s}_1^{-}
\end{equation}
where 
\begin{align}
  \mathcal{D}(\boldsymbol{\epsilon}^d) & := \{ (\boldsymbol{s}_1^{+},\boldsymbol{s}_1^{-},\boldsymbol{p}_1): \boldsymbol{s}_1^{+},\boldsymbol{s}_1^{-} \ge \boldsymbol{0}, \nonumber \\
  & \boldsymbol{D}\boldsymbol{x} + \boldsymbol{E}\boldsymbol{p}_1 + \boldsymbol{F}\boldsymbol{\theta}_1 + \boldsymbol{G}(\boldsymbol{\epsilon}^d + \boldsymbol{s}_1^{+} - \boldsymbol{s}_1^{-}) \le \boldsymbol{w} \label{eq_SPD_1} \\
  & \boldsymbol{e}^{\top}(\boldsymbol{p}^{g}_1 - \boldsymbol{p}^g) \le c^d \label{eq_SPD_3} \}.
\end{align}

The subproblem related to HLRU (SPR) is formulated as below with subscription 2 indicating associated variables.
\begin{equation}
  \text{(SPR)}\quad\max_{\boldsymbol{\epsilon}^r\in\mathcal{U}^r}\min_{(\boldsymbol{s}_2^{+},\boldsymbol{s}_2^{-},\boldsymbol{p}_2)\in\mathcal{R}(\boldsymbol{\epsilon}^r)}\boldsymbol{1}^{\top}\boldsymbol{s}_2^{+} + \boldsymbol{1}^{\top}\boldsymbol{s}_2^{-}
\end{equation}
where
\begin{align}
  \mathcal{R}(\boldsymbol{\epsilon}^r) & := \{ (\boldsymbol{s}_2^{+},\boldsymbol{s}_2^{-},\boldsymbol{p}_2): \boldsymbol{s}_2^{+},\boldsymbol{s}_2^{-} \ge \boldsymbol{0}, \nonumber \\
  & \boldsymbol{D}\boldsymbol{x} + \boldsymbol{E}\boldsymbol{p}_2 + \boldsymbol{F}\boldsymbol{\theta}_2 + \boldsymbol{G}(\boldsymbol{\epsilon}^r + \boldsymbol{s}_2^{+} - \boldsymbol{s}_2^{-}) \le \boldsymbol{w} \\
  & \boldsymbol{S}\boldsymbol{p}^g + \boldsymbol{T}\boldsymbol{p}^{g}_2 \le \boldsymbol{\Delta} \}. \label{eq_spr_redispatch}
\end{align}

The objectives of the above subproblems are to minimize the summation of non-negative slack variables, which indicate the un-accommodated uncertainty.

\subsection{Subproblem Reformulation} \label{sec_dual}
The above two subproblems are difficult to solve since there are infinite values of uncertain variables and re-dispatch strategies. However, it has been proved that the optimal solution is achieved at the extreme point of polygonal uncertainty set \cite{Ben-Tal2004,White1992}. Therefore, the subproblems are first converted into maximization problems according to the duality theory, and are then reformulated as MILP problems \cite{Ye2015c} by applying the closed form of extreme points.

\subsection{Relax-and-Enforce Decomposition (RED)}
Although the infinite values of continuous uncertain variables have been reduced to finite extreme points and the subproblems have been converted into MILP problems, the number of extreme points in large systems is so huge that the problems may be computationally intractable. Therefore, the RED technique is applied to decompose the two subproblems into smaller time-decoupled problems which can be solved efficiently. The only time-coupled constraint in SPD (\ref{eq_SPD_3}) is relaxed in the following relaxed subproblem SPD-1.
\begin{equation}
  \text{(SPD-1)}\quad\max_{\boldsymbol{\epsilon}^d\in\mathcal{U}^d}\min_{(\boldsymbol{s}_1^{+},\boldsymbol{s}_1^{-},\boldsymbol{p}_1)\in\mathcal{L}_1(\boldsymbol{\epsilon}^d)}\boldsymbol{1}^{\top}\boldsymbol{s}_1^{+} + \boldsymbol{1}^{\top}\boldsymbol{s}_1^{-} \nonumber
\end{equation}
where 
\begin{align}
  \mathcal{D}_1(\boldsymbol{\epsilon}^d) & := \{ (\boldsymbol{s}_1^{+},\boldsymbol{s}_1^{-},\boldsymbol{p}_1): \boldsymbol{s}_1^{+},\boldsymbol{s}_1^{-} \ge \boldsymbol{0}, \nonumber \\
  & \boldsymbol{D}\boldsymbol{x} + \boldsymbol{E}\boldsymbol{p}_1 + \boldsymbol{F}\boldsymbol{\theta}_1 + \boldsymbol{G}(\boldsymbol{\epsilon}^d + \boldsymbol{s}_1^{+} - \boldsymbol{s}_1^{-}) \le \boldsymbol{w} \}. \nonumber
\end{align}
Then it is enforced in SPD-2 as below.
\begin{equation}
  \text{(SPD-2)}\quad\max_{\boldsymbol{\epsilon}^d\in\mathcal{U}^l}\min_{\boldsymbol{p}_1\in\mathcal{L}_2(\boldsymbol{\epsilon}^d)} \boldsymbol{e}^{\top}(\boldsymbol{p}^{g}_1 - \boldsymbol{p}^g) - c^d \nonumber
\end{equation}
where 
\begin{align}
  \mathcal{D}_2(\boldsymbol{\epsilon}^d) := \{ \boldsymbol{p}_1: \boldsymbol{D}\boldsymbol{x} + \boldsymbol{E}\boldsymbol{p}_1 + \boldsymbol{F}\boldsymbol{\theta}_1 + \boldsymbol{G}\boldsymbol{\epsilon}^d \le \boldsymbol{w} \}. \nonumber
\end{align}

The SPD-1 and SPD-2 can both be formulated as MILP problems using the method in Section \ref{sec_dual}. Furthermore, since the annual net load duration curve does not represent the temporal relationship between different net load levels, they can be decomposed into subproblems for individual time slots and solved in parallel. 
Since the SPR is based on the annual net load duration curve, and there is no recourse cost constraint in it, it can simply be decomposed into time-independent subproblems. Meanwhile, the temporal relationship of net load ramping events in each subproblem is characterized using HLRU sets, while the re-dispatch limit is enforced by (\ref{eq_spr_redispatch}).

\begin{figure}[htb]
\centering
\includegraphics[width=3.5 in]{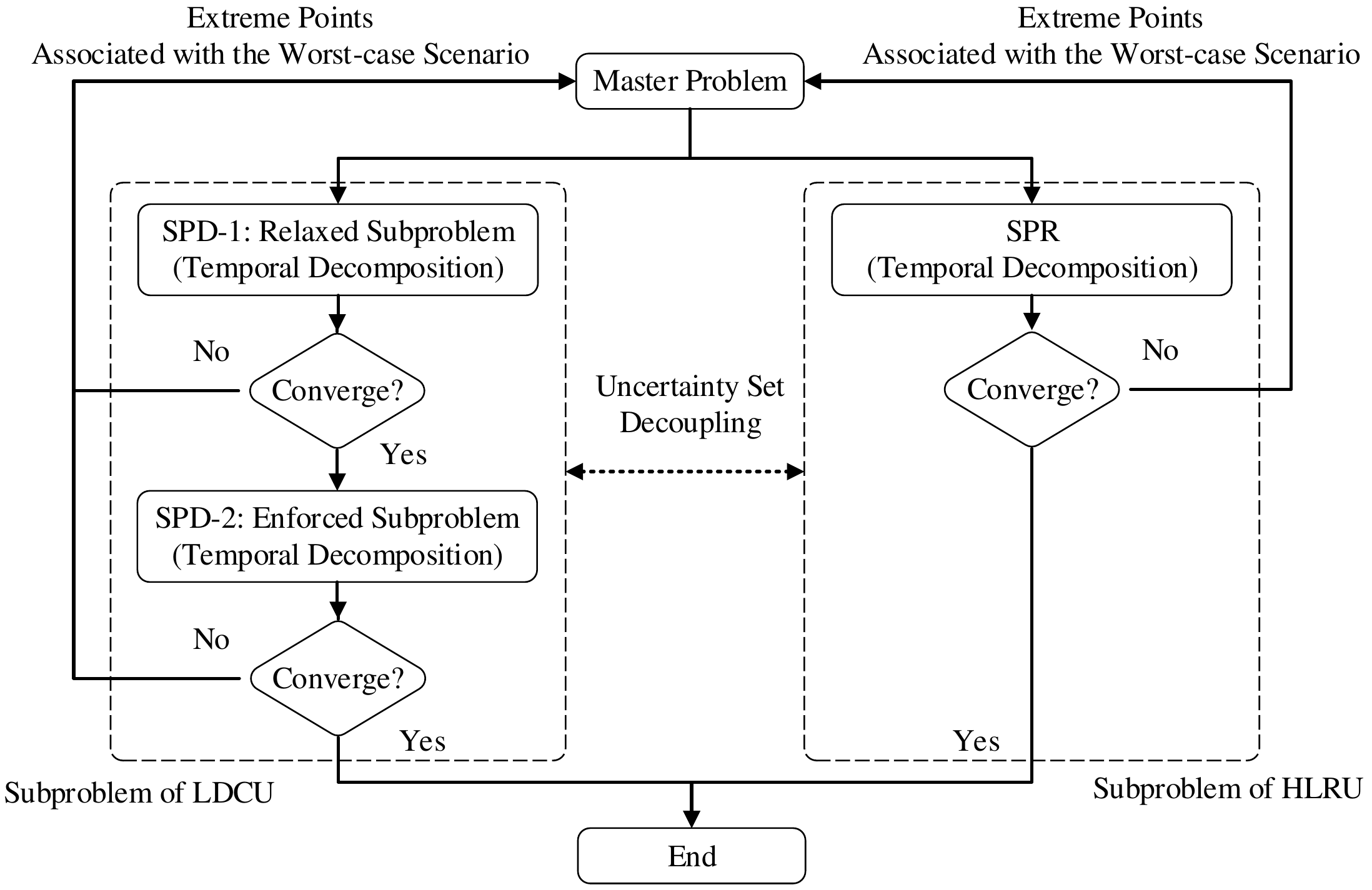}
\caption{Flowchart of the solution process.}
\label{fig_solution}
\end{figure}

The solution process is illustrated in Fig. \ref{fig_solution}, which contains one master problem and three subproblems, all formulated as MILP problems. The subproblems are solved sequentially and generate extreme points of uncertainty sets in each iteration. If the difference of the objective values of two iterations is within the pre-defined tolerance, the procedure is regarded to be converged. The convergence tolerance used in this paper is $10^{-3}$, which has been used in \cite{Ye2015c,Bertsimas2013}.

\section{Case Studies} \label{sec_case}

\subsection{A Modified Garver's 6-bus System}
\begin{figure}[htb]
\centering
\includegraphics[width=3 in]{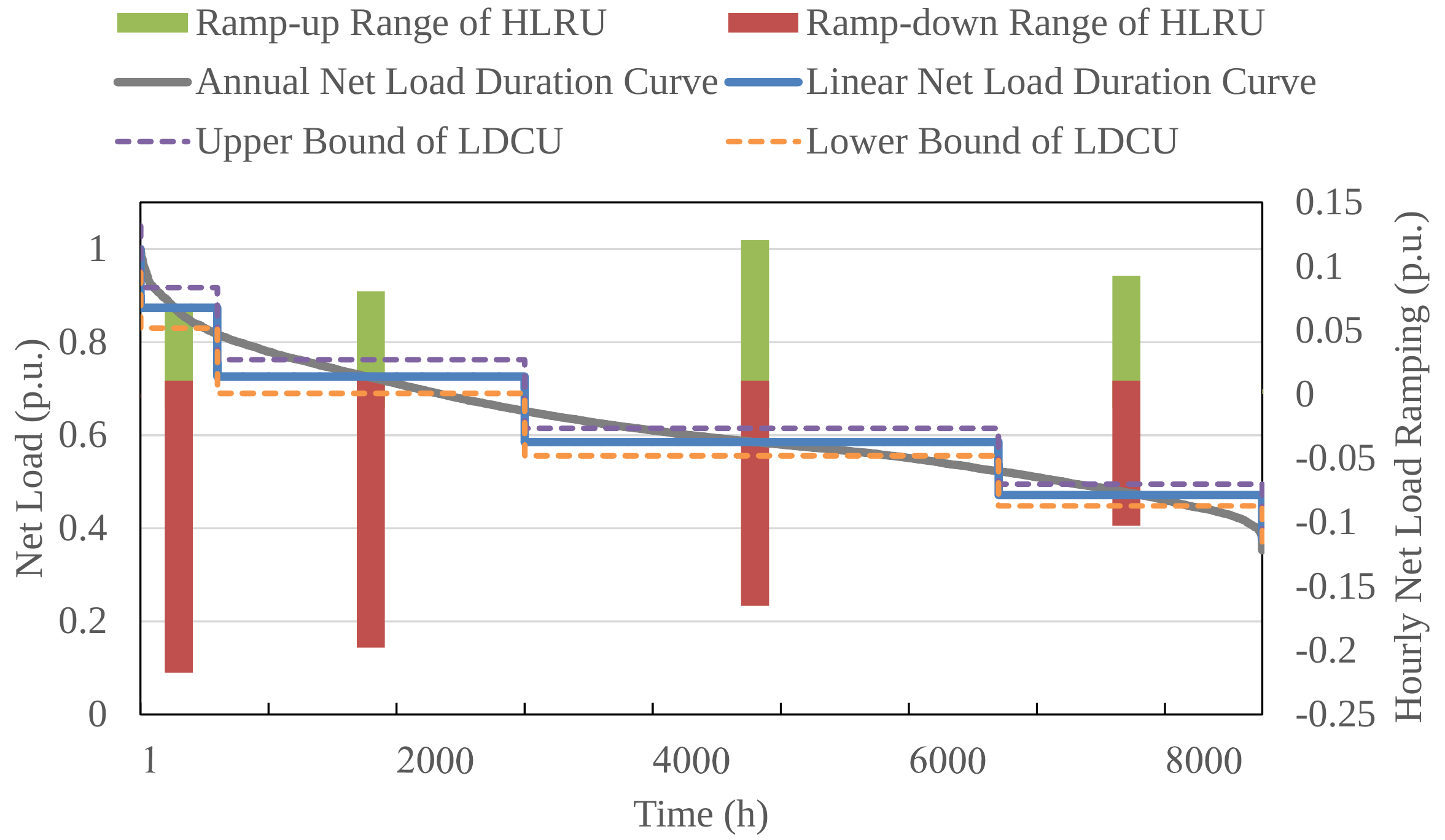}
\caption{Annual net load duration curve with two types of uncertainty.}
\label{fig_load}
\end{figure}
The proposed approach is applied to a modified Garver's 6-bus system \cite{Garver1970} where a line is added in corridor 2-6. The line data can be found in \cite{Alguacil2003a}. At most two parallel lines are allowed in each corridor. New generators can only be installed at bus 1, 3, 6, and at most two generators are allowed to install at each bus. Other detailed data can be found online at \url{http://motor.ece.iit.edu/Data/rotep}.

Real-world net load data from PJM in year 2015 \cite{PJM} are employed, as illustrated in Fig. \ref{fig_load}. The planning horizon is five years. The annual growth rate of net load and hourly net load ramping range are both assumed to be 5\%. The LDCU is assumed to be 5\% of each nodal net load.

\subsubsection{Flexibility of Different Resources}
The proposed planning model consists of three flexible resources: transmission lines (T), generators (G), and FACTS devices (F). Different combinations of these resources are compared to reveal some insights on how to coordinate them. The construction period of a transmission line is 1 year, while new generators and FACTS devices are assumed to be put into use in the same year when the investments are made. The coordination of construction period is discussed in Section \ref{sec_cp}. 

\begin{table}[htbp]
  \centering
\begin{threeparttable}

  \caption{Costs of Different Models}
    \begin{tabular}{lccccccc}
    \hline
    Cost (M\$)& T     & G     & F     & T+G\tnote{a}     & T+F   & G+F   & T+G+F \\
    \hline
    Line      & -\tnote{b}     & -     & -     & -   & -     & -     & 48.14 \\
    FACTS     & -     & -     & -     & -       & -     & -     & 15.46 \\
    Generator & -     & -      & -     & -   & -     & -     & 70.00 \\
    Investment\tnote{c} & -     & -     & -     & -  & -     & -     & 133.60 \\
    Operation\tnote{d} & -     & -      & -     & -  & -     & -     & 167.21 \\
    Total     & -     & -     & -     & -  & -     & -     & 300.81 \\
    \hline
    \end{tabular} \label{tab_cost}
    \begin{tablenotes}
        \item [a] The plus sign means the combination of different resources.
        \item [b] The hyphen indicates no result is obtained due to infeasibility.
        \item [c] Total investment cost.
        \item [d] Operation cost of base-case scenario.
    \end{tablenotes}
\end{threeparttable}   
\end{table}

As shown in Table \ref{tab_cost}, only T+G+F results in a feasible solution. Accordingly, we have the following observations from this case:

a) The master problems of models T, G, and F are infeasible. This implies that applying only one type of flexible resource cannot provide sufficient flexibility in this case.  

b) T+F is infeasible due to lack of new generators to meet the increasing load. Thus the role of generation-side flexibility may not be replaced by transmission flexibility.

c) G+F is infeasible due to lack of transmission capacity. FACTS devices usually have smaller capacity than lines and their installation relies on lines. Therefore FACTS devices may not be completely take the place of transmission lines in planning, but rather serve as a supplement.

d) T+G is infeasible since the LDCU cannot be coped with in this case, showing that more flexibility is required to handle LDCU.

e) Only T+G+F is capable of addressing LDCU and HLRU simultaneously, implying that LDCU and HLRU lead to much more requirement of flexibility. 

These observations evidently reveal that it is necessary to coordinate different types of flexible resources in order to accommodate LDCU and HLRU at the same time.

\subsubsection{Coordination of Construction Periods}\label{sec_cp}
In practice, the construction of transmission lines and generators usually takes more time than that of wind farms, leading to mismatch in expansion. In this section, model T+G+F is tested with different construction periods of transmission lines and generators, to show the impact of construction periods on planning strategy and the underlying benefits of FACTS.

Table \ref{tab_cost_ct} shows the costs with different construction periods of lines, while that of generators and FACTS are 0 year. When the construction period of lines increases from 0 year to 1 year, the investment of line decreases, while the investment of FACTS increases, indicating that FACTS devices are capable of providing transmission flexibility to coordinate with longer-term construction of line. Note that higher operation cost of base-case scenario and total cost are also experienced, implying that FACTS devices are unable to completely replace transmission lines. In particular, when the construction period of line becomes 2 years, the problem is infeasible as the LDCU cannot be accommodated fully. 

\begin{table}[htbp]
  \centering
  \caption{Costs with Different Construction Period of Line}
    \begin{tabular}{lccc}
    \hline
          & \multicolumn{3}{c}{Construction Period (Year)} \\
    \hline
    \multicolumn{1}{l}{Cost (M\$)} & 0     & 1     & 2 \\
    \hline
    \multicolumn{1}{l}{Line} & 60.00 & 48.14 & - \\
    \multicolumn{1}{l}{FACTS} & 0.00  & 15.46 & - \\
    \multicolumn{1}{l}{Generator} & 70.00 & 70.00 & - \\
    \multicolumn{1}{l}{Investment} & 130.00 & 133.60 & - \\
    \multicolumn{1}{l}{Operation} & 150.41 & 167.21 & - \\
    \multicolumn{1}{l}{Total} & 280.41 & 300.81 & - \\
    \hline
    \end{tabular}
  \label{tab_cost_ct}
\end{table}

When the construction period of generators increases to 1 year, the problem is infeasible even when the construction period of lines is 0 year, since no new generator can be installed in the first year to address the HLRU. It demonstrates the irreplaceable role of generation-side flexibility in the planning problem with uncertainty.

Overall, it may be understood that FACTS devices are supplemental transmission flexible resources that can help coordinate construction periods, while the generation-side flexibility has a crucial role in addressing the HLRU. Besides, shorter construction periods of conventional resources can better facilitate the integration of wind power.

\subsubsection{Consideration of Uncertainty}\label{subsubsection_uncertainty}
In this section, the following three models are compared to analyze the impact of uncertainty. Note that M1 and M2 are indeed parts of M3. The costs of these models are listed in Table \ref{tab_cost_uncertainty}.

M1: the planning model not considering uncertainty.

M2: the planning model considering LDCU.

M3: the proposed model considering LDCU and HLRU.

\begin{table}[htbp]
  \centering
  \caption{Costs Under Different Uncertainties}
    \begin{tabular}{lccc}
    \hline
    Cost (M\$) & M1    & M2    & M3 \\
    \hline
    Line  & 60.00 & 60.00 & 48.14 \\
    FACTS & 4.00  & 9.00 & 15.46 \\
    Generator & 25.00 & 50.00 & 70.00 \\
    Investment & 89.00 & 119.00 & 133.60 \\
    Operation & 146.78 & 147.75 & 167.21 \\
    Total & 235.78 & 266.75 & 300.81 \\
    \hline
    \end{tabular}
  \label{tab_cost_uncertainty}
\end{table}

It is found that the investment cost, the operation cost of base-case scenario, and the total cost increase as more uncertainty is considered. The details of planning strategy are given in Table \ref{tab_strategy_case6}. Accordingly, we have the following observations:

a) \emph{Lines:} The planning strategies of the three models all consist of at least one new line in corridor 4-6, showing the importance of the investment in line.

b) \emph{FACTS:} Compared with M1, in M2 where only LDCU is considered, FACTS devices with larger capacity (type II) is installed. Moreover, in M3 where HLRU is also taken into account, not only the total FACTS capacity increases, but more FACTS devices are installed in the 4th year, leading to a postpone of new line construction. It demonstrates that FACTS devices considerably help provide transmission flexibility, especially when HLRU is considered.

c) \emph{Generators:} When considering HLRU, a new generator with larger ramping capability (type II) is installed, showing there is additional requirement on generation-side flexibility.

\begin{table}[htbp]
  \centering
\begin{threeparttable}
  \caption{Planning Strategies with Different Models}
    \begin{tabular}{clccccc}
    \hline
          & Year  & 1     & 2     & 3     & 4     & 5 \\
    \hline
    \multirow{3}{*}{M1} & Line  & -\tnote{a}      & 4-6 4-6\tnote{b}  & -      & -     & -  \\
          & FACTS & 1-2(I) 1-4(I)\tnote{c} & -      & -      & -      & - \\
          & Generator & 1\tnote{d} (I) & -       &-       & -      & -  \\
    \hline
    \multirow{4}{*}{M2} & Line  & -      & 4-6 4-6  & -      & -       & -  \\
          & \multirow{2}{*}{FACTS} & 1-2(I) 1-4(II)  & -  & -        & -      &-  \\
          & & 2-4(I) & -  & -        & -      &-  \\
          & Generator & 1(I) 3(I) & -      &  -     &-       & - \\
    \hline
    \multirow{5}{*}{M3} & Line  & -      & 4-6   &  -     & 1-5      &-  \\
          & \multirow{2}{*}{FACTS} & 1-2(I) 1-4(II)  & -  & -        & 2-6(I) 4-6(I)      &-  \\
          & & 2-4(II) & -  & -        & -      &-  \\
          & Generator & 1(II)  & -      &  -     & -     & - \\
    \hline
    \end{tabular} \label{tab_strategy_case6}
    \begin{tablenotes}
        \item [a] The hyphen indicates no new facility is built or installed.
        \item [b] Corridor number, each of which indicates a new line/FACTS device.
        \item [c] The number in parenthesis indicates the type of FACTS/generator.
        \item [d] Bus number, each of which indicates a new generator.
    \end{tablenotes}
   \end{threeparttable}
\end{table}

In order to investigate the impacts of different uncertainties on operational feasibility, we test the planning strategy of M1 with feasibility sets $\mathcal{F}^d$ and $\mathcal{F}^r$. In the test of $\mathcal{F}^d$, the planning strategy is feasible when the transmission capacity constraints on lines 1-4 and 1-5 are relaxed, and the capacity constraint of either G1 or G3 is relaxed. As the HLRU is considered at the same time, besides the aforementioned line and generator capacity constraints, the ramping constraint of at least one of the existing generators has to be relaxed. In short, the consideration of LDCU and HLRU leads to more requirement of transmission and generation-side flexibility in this case.

As to M3, when the convergence tolerance is set to $10^{-3}$, 3 cutting planes are generated by SPD, while 2 cutting planes are generated by SPR. When the convergence tolerance decreases to $10^{-4}$, the number of cutting planes generated by SPD increases rapidly to 10, while that generated by SPR remains unchanged. However, the objective value remains unchanged. When the convergence tolerance increases to $10^{-2}$, the number of cutting planes is the same as that with the tolerance of $10^{-3}$, and the objective value is unchanged. Therefore, the convergence tolerance of $10^{-3}$ is accurate enough for this case.

\begin{figure}[htb]
\centering
\includegraphics[width=3.2in]{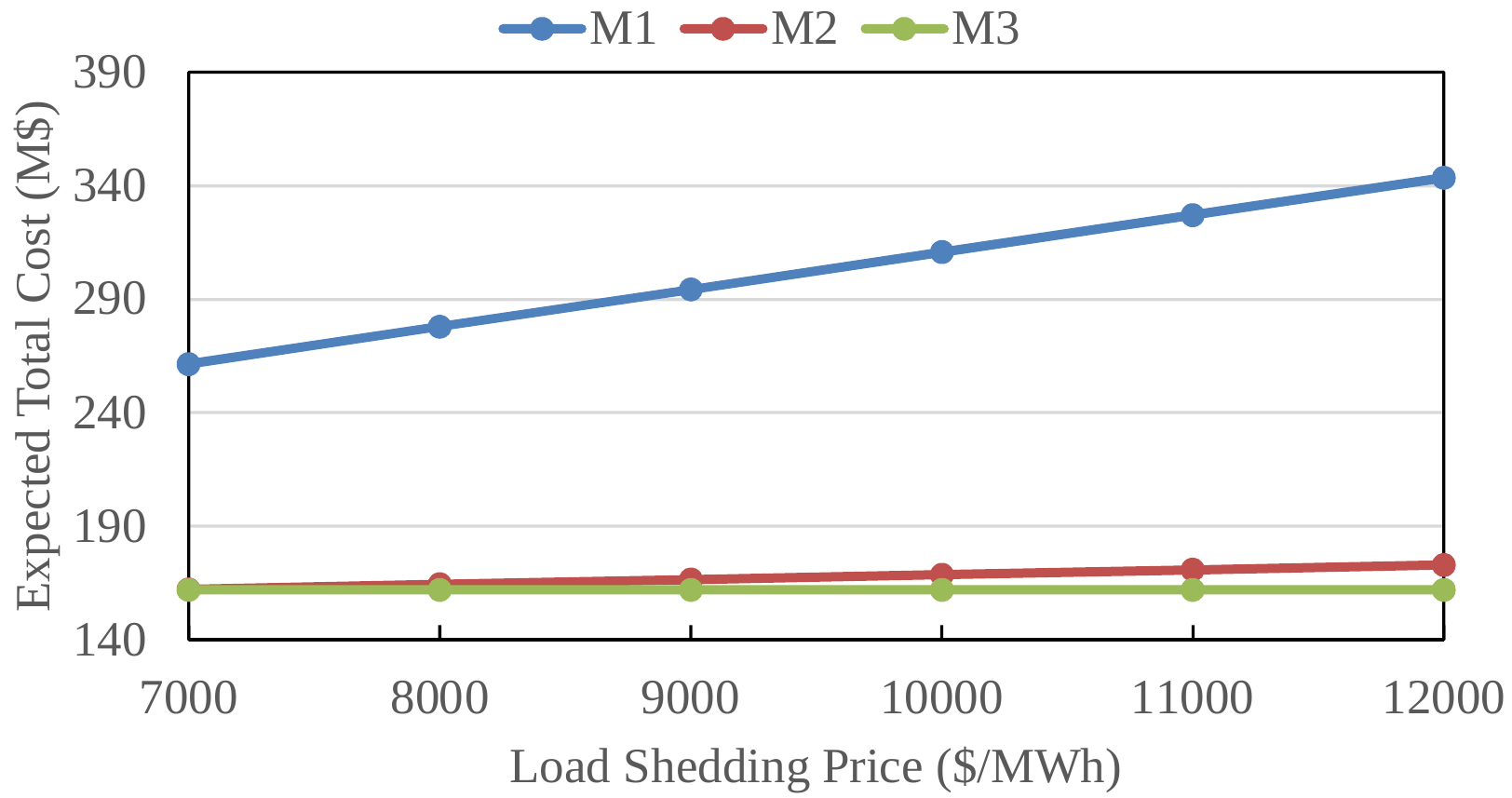}
\caption{Expected total cost with different load shedding prices.}
\label{fig_lsp}
\end{figure}

To test the robustness of the solutions of planning strategy and unit commitment schedule, dispatch simulations are performed for two sets of 1000 randomly generated scenarios. One set corresponds to the LDCU, which follows a normal distribution with the mean $P^d_{i,y,h}$ and the standard deviation $\bar{u}_{i,y,h}^d$, and is assumed to be within the interval $[P^d_{i,y,h}-\bar{u}_{i,y,h}^d, P^d_{i,y,h}+\bar{u}_{i,y,h}^d]$. The other set is related to the HLRU, which follows a uniform distribution in the interval $[P^d_{i,y,h}+\underline{u}_{i,y,h}^r, P^d_{i,y,h}+\bar{u}_{i,y,h}^r]$. Load shedding is allowed at the price of \$7,000/MWh. When the LDCU and HLRU scenarios are considered at the same time, the power demand in both scenarios are met, while the ramping constraints between the two scenarios are satisfied. Simulation results are provided in Table \ref{tab_sim}. It can be observed that when only the LDCU is considered, M2 achieves the least ETC. The ETC of M3 is higher than M2 since it is more conservative. When the LDCU and HLRU are both taken into account, though the EOC of M3 is higher than other models since only the operation cost of the base-case scenario is minimized, M3 outperforms the other models in terms of ETC. In particular, M3 avoids the extremely high load shedding cost in the possible scenario, indicating the necessity to consider the HLRU. Since the average load shedding price for a developed, industrial economy range from approximately \$9,000/MWh to \$45,000/MWh \cite{LondonEconomicsInternationalLLC2013}, M3 has a bigger advantage when the load shedding price is higher, as shown in Fig. \ref{fig_lsp}. In short, the robust model has a great advantage over the deterministic model when the load shedding price is high, achieving lower ETC and avoiding the high cost in the worst-case scenario.

\begin{table}[htbp]
  \centering
  \begin{threeparttable}
  \caption{Simulation Results with LDCU and HLRU}
    \begin{tabular}{ccccccccc}
    \hline
    & \multirow{2}{*}{Model} & ETC\tnote{a} & EOC & ELC & HLC\tnote{b} & EENS\tnote{c} \\
    &    & (M\$) & (M\$) & (M\$) & (M\$) & MWh/year\\
    \hline
    \multirow{3}{*}{LDCU} & M1 & 176.00  & 146.78  & 29.21 & 216.58 & 1014.36\\
    & M2    & 147.26 & 147.26 & 0 & 0 & 0 \\
    & M3    & 162.03  & 162.03 & 0  & 0 & 0 \\
    \hline
    LDCU & M1 & 261.53  & 146.80  & 114.73 & 518.51 & 3936.96 \\
    + & M2    & 162.30 & 147.28 & 15.02 & 472.86 & 488.67 \\
    HLRU & M3    & 162.04 & 162.04 & 0 & 0 & 0 \\
    \hline
    \end{tabular}
  \label{tab_sim}
  \begin{tablenotes}
        \item [a] Expected total cost, including expected operation cost (EOC) and expected load shedding cost (ELC).
        \item [b] Highest load shedding cost among all the scenarios.
        \item [c] Expected energy not supplied.
    \end{tablenotes}
  \end{threeparttable}
\end{table}

Since it is difficult to forecast hourly load accurately several years ahead, we randomly generate 10 scenarios each year based on the hourly load data of PJM from year 2011 to 2015 (50 scenarios in total). Each scenario contains 8760 hours load data, and the scenarios are assumed to follow normal distributions. Daily unit commitment and economic dispatch are performed for these scenarios.

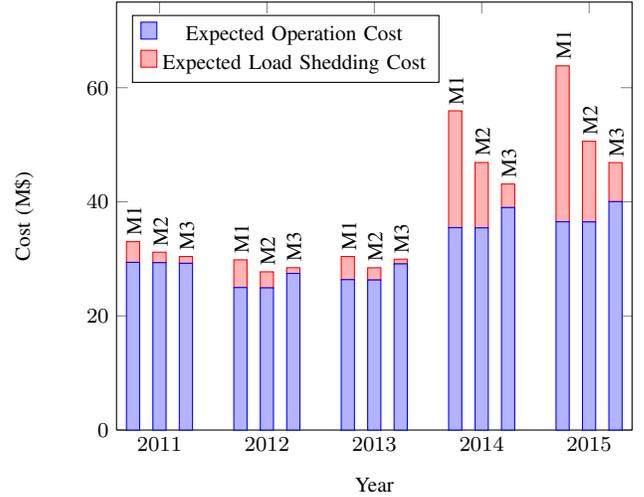
\begin{figure}[htb]
\centering
\begin{tikzpicture}
\begin{axis}[nodes near coords ybar stacked configuration/.style={},
  ybar stacked,
  x tick label style={
    /pgf/number format/1000 sep=},
  bar shift=-10pt,
  bar width=5pt,
  ymin=0,
  ymax=75,
  ylabel=Cost (M\$), 
  xlabel=Year,
  legend style={legend pos=north west},
  ]
\addplot coordinates
{(2011,29.38) (2012,24.97) (2013,26.36) (2014,35.47) (2015,36.50)};
\addplot+[
  nodes near coords=M1,
  every node near coord/.append style={
    black,
    xshift=-4pt,
    yshift=8pt,
    rotate=90,
    font=\footnotesize},
] coordinates
{(2011,3.64) (2012,4.86) (2013,4.05) (2014,20.48) (2015,27.33)};
\legend{Expected Operation Cost, Expected Load Shedding Cost}
\footnotesize
\end{axis}

\begin{axis}[nodes near coords ybar stacked configuration/.style={},
  ybar stacked,
  bar shift=0pt,
  bar width=5pt,
  ymin=0,
  ymax=75,
  axis lines=none,
  ]  
\addplot coordinates
{(2011,29.34) (2012,24.92) (2013,26.30) (2014,35.43) (2015,36.48)};
\addplot+[
  nodes near coords=M2,
  every node near coord/.append style={
    black,
    xshift=6pt,
    yshift=8pt,
    rotate=90,
    font=\footnotesize},
] coordinates
{(2011,1.83) (2012,2.81) (2013,2.14) (2014,11.46) (2015,14.13)};
\footnotesize
\end{axis}

\begin{axis}[nodes near coords ybar stacked configuration/.style={},
  ybar stacked,
  bar shift=10pt,
  bar width=5pt,
  ymin=0,
  ymax=75,
  axis lines=none,
  ]  
\addplot coordinates
{(2011,29.22) (2012,27.45) (2013,29.12) (2014,38.99) (2015,40.03)};
\addplot+[
  nodes near coords=M3,
  every node near coord/.append style={
    black,
    xshift=16pt,
    yshift=8pt,
    rotate=90,
    font=\footnotesize},
] coordinates
{(2011,1.19) (2012,1.00) (2013,0.83) (2014,4.13) (2015,6.84)};
\footnotesize
\end{axis}
\end{tikzpicture}
\caption{Expected costs of different models.}
\label{fig_cost3}
\end{figure}

The ETC, including the EOC and the ELC, is provided in Fig. \ref{fig_cost3}. The results generally coincide with our previous simulations. The solution of M3 ensures robustness to a large extent, though there is still load shedding in M3, mainly because in the optimization model, the uncertainty sets are constructed based on the stepwise net load duration curves and the ramping events are assumed to start from each average net load level of the stepwise net load duration curves.

The deterministic model, M1, results in the highest ETC since no uncertainty is considered. Note that the difference between M2 and M3 is that M2 only considers LDCU, while M3 considers both LDCU and HLRU. Comparing M3 with M2, it is observed that the ETC of M3 is lower than that of M2 in most of the years, especially in years 2014 and 2015. The average ETC in 5 years is reduced by about 3.3\%. Besides, the ELC of M3 is apparently lower than that of M2, showing the effectiveness of M3 in capturing the ramping uncertainty. The advantage of M3 in reducing ELC is especially evident in years 2014 and 2015 when the ELC is higher than that in other years. 

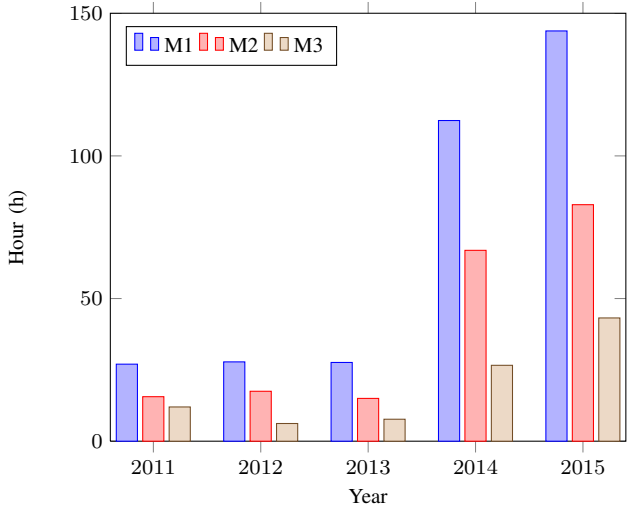
\begin{figure}[htb]
\centering
\begin{tikzpicture}
\begin{axis}[
  ybar,
  ymin=0,ymax=150,
  x tick label style={
    /pgf/number format/1000 sep=},
  bar width=8pt,
  legend style={legend pos=north west,legend columns=-1},
  xlabel={Year},
  ylabel={Hour (h)},
]
  \addplot coordinates {(2011,27) (2012,27.8) (2013,27.6) (2014,112.4) (2015,143.8)};
  \addplot coordinates {(2011,15.6) (2012,17.5) (2013,15) (2014,66.9) (2015,82.9)};
  \addplot coordinates {(2011,12) (2012,6.2) (2013,7.7) (2014,26.6) (2015,43.2)};
  \legend{M1, M2, M3}
  \footnotesize
\end{axis}
\end{tikzpicture}
\caption{Loss of load hours of different models.}
\label{fig_lolh3}
\end{figure}

\begin{figure}[htb]
\centering
\begin{tikzpicture}
\begin{axis}[
  ybar,
  ymin=0,
  ymax=3000,
  bar width=8pt,
  x tick label style={
    /pgf/number format/1000 sep=},
  legend style={legend pos=north west,legend columns=-1},
  xlabel={Year},
  ylabel={Energy (MWh)},
]
  \addplot coordinates {(2011,364.35) (2012,485.95) (2013,404.72) (2014,2047.81) (2015,2733.04)};
  \addplot coordinates {(2011,182.82) (2012,280.82) (2013,213.93) (2014,1146.11) (2015,1412.90)};
  \addplot coordinates {(2011,118.78) (2012,100.39) (2013,82.81) (2014,412.58) (2015,684.33)};
  \legend{M1, M2, M3}
  \footnotesize
\end{axis}
\end{tikzpicture}
\caption{Expected energy not supplied of different models.}
\label{fig_eens3}
\end{figure}
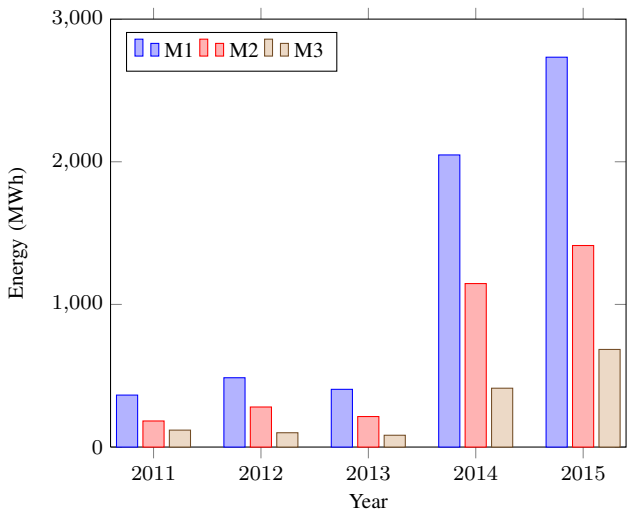

In order to further investigate the effect of M3, the loss of load hours (LOLH) and EENS are provided in Fig. \ref{fig_lolh3} and \ref{fig_eens3}, respectively. The LOLH of M3 is about half of that of M2 every year. Similar situations can be observed with respect to the EENS. It implies that when the uncertainty of ramping is high, causing a lot of load shedding, M3 outperforms the other models apparently, achieving less expected load shedding and a lower expected total cost.

\subsection{IEEE 118-bus System}
To evaluate the practicality of the proposed model and solution method, we further conduct tests on the IEEE 118-bus system with 10 selected corridors to build new lines, 5 selected corridors to install FACTS devices, and 5 selected buses to install new generators. The peak total net load in the first year of a 5-year planning horizon is 6500MW. Other detailed data can be found at \url{http://motor.ece.iit.edu/Data/rotep}. The proposed models are implemented in GAMS \cite{GAMS} and solved using CPLEX \cite{CPLEX}. Table \ref{tab_118} provides the results with different number of buses associated with uncertainty. It can be observed that the total cost increases as more uncertainty is considered, and the proposed solution approach (Alg.1) improves computational efficiency over the C\&CG approach without the RED technique (Alg.2), where the max-min subproblems are formulated as MILP problems using the method described in Section \ref{sec_dual} and then solved with CPLEX.

\begin{table}[htbp]
  \centering
  \begin{threeparttable}
  \caption{Objective Value and Computational Time of Different Models}
    \begin{tabular}{ccccccc}
    \hline
          & No. of & \multicolumn{2}{c}{Objective Cost (M\$)}  & \multicolumn{2}{c}{CPU Time (s)} \\
          & Buses\tnote{a} & Alg.1 & Alg.2 & Alg.1 & Alg.2 \\
    \hline
    M1    & -      & \multicolumn{2}{c}{3001.23} & \multicolumn{2}{c}{56.27} \\
    \hline
    \multirow{4}{*}{M2}    & 5     & 3137.65 & 3138.68 & 477.64 & 466.34 \\
         & 10     & 3137.91 & 3138.80 & 200.86 & 2527.56 \\
          & 15     & 3146.88 & 3147.06 & 379.43 & 702.10 \\
          & 20     & 3148.73 & 3147.37 & 388.75 & 2475.01 \\
    \hline
    \multirow{4}{*}{M3}    & 5      & 3138.56 & 3138.54 & 755.50 & 1159.13 \\
          & 10     & 3165.87 & 3164.01 & 456.19 & 5758.36 \\
          & 15     & 3235.08 & 3234.29 & 1221.83 & 4081.50 \\
         & 20     & 3340.64 & 3340.84 & 1683.01 & 13662.31 \\
    \hline
    \end{tabular}%
  \label{tab_118}%
  \begin{tablenotes}
        \item [a] {Number of buses with uncertainty.}
    \end{tablenotes}
  \end{threeparttable}
\end{table}

Note that the time consumption of the master problem becomes larger during iterations, as more and more extreme points are added into it, especially in the iterations of SPR where the extreme points generated from SPD are also included. It should be noted that there is no strict requirement of computational time in planning problems. Moreover, the SPD and the SPR can be solved in parallel to accelerate the solution process, which are solved sequentially in this paper.

\subsection{Real System of Gansu Province in China}
To further evaluate the scalability of the proposed methodology, tests based on the data of the real system of Gansu province in China are conducted. The system contains 157 buses, 258 lines and 20 wind farms. The tests are conducted with 10 selected corridors to build new lines, 6 selected corridors to install FACTS devices, and 5 selected buses to install new generators. The planning horizon is 5 years. Table \ref{tab_gs157} provides the results.
\begin{table}[htbp]
  \centering
  \caption{Objective Value and Computational Time of Different Models}
    \begin{tabular}{cccccc}
    \hline
          & \multicolumn{2}{c}{Objective Cost (M\$)}  & \multicolumn{2}{c}{CPU Time (s)} \\
           & Alg.1 & Alg.2 & Alg.1 & Alg.2 \\
    \hline
    M1  & \multicolumn{2}{c}{11820.01} & \multicolumn{2}{c}{36.24} \\
    \hline
    M2  & 11936.99 & 11936.99 & 415.50 & 8392.53 \\
    \hline
    M3  & 12061.76 & 12061.78 & 4225.07 & 24253.62 \\
    \hline
    \end{tabular}
  \label{tab_gs157}
\end{table}

It is observed that Alg.1 is still faster than Alg.2 when applied in the real-world system. The advantage is as much as one order of magnitude. However, the difficulties of applying the proposed methodology in solving real systems may lie in the large scale of the systems and the large number of candidate sites for new devices, leading to high computational burdens on the iteration process.

\section{Concluding Remarks} \label{sec_conclusion}
The ever increasing uncertain and volatile wind power generation enforce new requirement of ramping flexibility and coordination of construction periods in expansion planning. To address such issues, we have developed a comprehensive robust planning model incorporating different flexible resources with different construction periods. We have constructed a novel uncertainty set to depict the HLRU in addition to the LDCU. The two types of uncertainty are decoupled and solved by combining the C\&CG and the RED algorithms. Real-world data is used to verify the effectiveness of the proposed model. Some remarks are provided as below.

1) The comprehensive multi-year planning model has provided some insights for power system planners to coordinate different resources and construction periods under uncertainty. FACTS is an effective supplemental tool to help coordinate different construction periods, since it usually takes much less time to install and provides transmission flexibility, but it cannot fully take the place of transmission lines. Generators provide generation-side flexibility to balance load and hedge against load ramping.

2) The proposed novel uncertainty sets capture the ramping uncertainty to a large extent while keeping the optimization problem tractable. By taking the HLRU into account, a more robust decision can be made with an effective reduction of expected load shedding, while the possible extremely high load shedding cost in the worst-case scenario is avoided. The benefits become greater when the load shedding price is higher, implying that the proposed model is suitable for high reliability requirement.

3) Even though the proposed multi-year planning problem has been formulated as a tractable two-stage robust optimization problem, the computational burden is still high, since two coupled uncertainties are considered. Decoupling and decomposition techniques have significantly improved computational efficiency. Tests on the
IEEE 118-bus system and a real-world system have shown that the acceleration of the solution process can be
as much as one order of magnitude compared with the standard C\&CG method without a temporal decomposition.




In our future work, some reliability indexes will be introduced to control the conservativeness of the model. The reliability constraints can be added on the worst-case scenario of the subproblems, using the formulation similar to the recourse cost constraints. Besides, unit commitment will be performed in the subproblems to better simulate daily operation.

\ifCLASSOPTIONcaptionsoff
  \newpage
\fi

\bibliographystyle{IEEEtran}
\bibliography{library}

\end{document}